# Clear thinking, vague thinking and paradoxes

Arieh Lev and Gil Kaplan

## Abstract

Many undergraduate students of engineering and the exact sciences have difficulty with their mathematics courses due to insufficient proficiency in what we in this paper have termed **clear thinking**. We believe that this lack of proficiency is one of the primary causes underlying the common difficulties students face, leading to mistakes like the improper use of definitions and the improper phrasing of definitions, claims and proofs. We further argue that clear thinking is not a skill that is acquired easily and naturally – it must be consciously learned and developed.

The paper describes, using concrete examples, how the examination and analysis of classical paradoxes can be a fine tool for developing students' clear thinking. It also looks closely at the paradoxes themselves, and at the various solutions that have been proposed for them. We believe that the extensive literature on paradoxes has not always given clear thinking its due emphasis as an analytical tool. We therefore suggest that other disciplines could also benefit from drawing upon the strategies employed by mathematicians to describe and examine the foundations of the problems they encounter.

## Introduction

Students working towards a B.A. in the exact sciences and in engineering encounter significant difficulties in their mathematics courses from day one. These students must be able to read, understand and accurately repeat definitions, and to skillfully and rigorously read and write mathematical claims and laws. In other words, they must be skilled in what this paper will refer to as "clear thinking": thinking that requires one to fully understand what concepts mean, to be able to use those concepts accurately, and to understand and express definitions and claims fully and accurately. Students must meet the requirements of clear thinking throughout all the stages involved in addressing a problem: in reading and understanding it, in analyzing it, and in formulating a solution.

The difficulties students face arise from three primary causes:

1. Clear thinking is different from the thinking we have been accustomed to engage in since birth. We are used to employing intuitive-associative thinking, which is very powerful, but which does not include the necessary control, or the attempt to carefully and accurately explore the meanings of concepts and claims. As a result, this thinking is susceptible to bias and error (see for instance [EW], [K], [L]). Clear thinking does not come automatically to anyone; most people require extensive training and practice to achieve satisfactory skill levels ([L]).

2. In many cases, students begin their studies without having undergone any training in clear thinking. Mathematics education in schools often emphasizes calculation over thinking, and practicing the use of clear thinking is generally overlooked in other fields as well.

3. Very early in their academic studies, students in mathematics courses are required to work with definitions and proofs at high levels of difficulty and complexity. There is something of a hidden assumption in academic institutions, according to which the students who arrive there have the necessary clear thinking skills to handle the demands of the coursework. It appears that, for a high percentage of these students, this assumption is incorrect ([EW], [L]).

These three weighty issues are joined by an additional, "psychological" problem: students are often unsure of what the point and purpose of clear thinking is. After all, it may seem easier (and even more practical) to make do with an intuitive understanding of concepts (like limits and continuity of functions), instead of getting "bogged down" in complicated definitions that make topics that seem clear and intuitive seem complex and difficult to understand. Such an attitude on students' part is unsurprising, especially if we take into account the historical development of many basic mathematical concepts. Thus, for instance, for many years mathematicians calculated limits and derivatives in ways that would today be considered (and that were considered by some in the past as well) to be lacking in any firm and rigorous basis. Only later generations of mathematicians were able to formulate the rigorous definitions that we find acceptable today. And if this was the case for mathematicians, how could we expect otherwise of first year students?

One central goal of this article is to propose a method for developing clear thinking skills through the examination and analysis of paradoxes. This method seems appropriate to us for the following reason: many paradoxes are simple to present, but can nevertheless elicit confusion and cognitive dissonance. Moreover, they are easier to understand and solve once the concepts and definitions used by the paradox are clearly and accurately described. This means that paradoxes can be used:

1. To introduce students simply, visibly and compellingly to the advantages of clear thinking.

2. To train students in accurately formulating definitions and claims, and in distinguishing between such accurate formulations on the one hand and phrasing that is vague and ambiguous on the other. (In other words: to help students acquire clear thinking skills.)

3. To show how easily concepts' lack of clarity can draw us into a cycle of confusion and contradiction.

4. To achieve the goals listed above simply and directly, without any need to first introduce the students to complex concepts and structures.

It is important to note at this point that practice and experience are necessary prerequisites for acquiring these skills. It is therefore important that the students first be required to face the problems (i.e. the paradoxes) on their own (with proper guidance), so that the detailed discussion of the problems is conducted only after the students have read, worked, thought and tried to reach conclusions by themselves.

In addition to the goal stated above, this article also has the secondary goal of drawing attention to the importance of clear thinking – in all fields, not just in mathematics and the exact sciences. This is true, for instance, in the field of paradoxes itself. In reading the books and articles that discuss paradoxes, we found that many do not pay sufficient attention to clearly analyzing the claims and concepts that appear in the paradox, as a preliminary step to be carried out *before* continuing to discuss it.  We believe it would be better if those who addressed these issues borrowed some practices from the mathematicians and applied them to paradox analysis from its earliest stages. This would make it possible to remove some of the vagueness and confusion at the preliminary stage, making it easier to continue analyzing the problem later.

The primary source we will be using on the topic of paradoxes is the book *Paradoxes*, by R.M. Sainsbury [Sa]. We chose this book because in our opinion (and that of many others) it offers the clearest and most faithful review of the extensive literature available on the topic of paradoxes. That said, when on occasion we criticize a particular approach to a topic, our critique will refer to the approach itself, and not to any claim or manner of presentation on the part of the author who described it. For the purpose of our discussion, we quote the definition of a paradox that appears in [Sa, p. 1], namely:

**An apparently unacceptable conclusion derived by apparently acceptable reasoning from apparently acceptable premises.**

# A: Preliminary presentation: "deceptive" vs. "true" paradoxes

Let us begin with two ancient, famous paradoxes: "the Liar" (we will discuss only a basic, simple version of this) and "Achilles and the tortoise." These paradoxes will convincingly convey to the students how a preliminary discussion to clarify the concepts and claims makes it far easier to clarify problems that may seem daunting and confusing at first glance.

**The Liar:** In the simple version of this paradox, we must examine the words of a man who says of himself, "I am a liar," and ask: is this man telling the truth or is he lying? He cannot be telling the truth, since then his testimony that he is a liar is a lie. Therefore, he is a liar. But this also cannot be, since then his testimony that he is a liar would be truth, which means he cannot be a liar. We must therefore conclude that he is neither lying nor telling the truth, which is impossible! It appears that this simple version of the paradox leaves many people in a state of perplexity and confusion. But there is a simple way of avoiding this state: all we must do is clearly examine the meaning of the concepts in this paradox (i.e. "liar" and "telling the truth"). We do this by trying to formulate specific definitions for these concepts, like those presented below.

<u>Definition:</u>

a) We say that a man is a **truth teller** if everything he says is always true.

b) We say that a man is a **liar** if everything he says is always a lie.

Now all we must do is carefully test the claims we have raised in our analysis of the paradox, in light of these definitions. If we do so, we will reach the conclusion that our perplexity is unjustified: a man who is not a truth teller and not a liar is a man who neither tells the truth all the time, nor lies all the time. And this is most certainly possible! In fact, the paradox describes a person who is like all of us – sometimes truthful, sometimes lying (as he is when he makes the claim, "I am a liar") or even capable sometimes of making two contradicting claims. Therefore, there is no contradiction here, and no paradox, but rather a confusing presentation of the problem. Moreover, the discussion above shows us that the way out of this confusion is through clear examination of the concepts and the claims.

Some might counter by saying: what if this refers to a (hypothetical) world in which every person is either a truth teller or a liar in the senses defined above (as is the case in many popular riddles)? Will it be a paradox then? Before answering this query, we must note that if this was indeed the intention, then the problem - as written above - is invalid. It is improper to present a problem that is based upon a hidden assumption. All assumptions must be stated explicitly in full within the problem itself! Now let us move on to the answer. The case raised by these students posits a model in which the following basic assumption (axiom) is true:

Assumption: Every person in this model is either a truth teller or a liar.

Now (assuming that this model is indeed consistent) we would ostensibly get a paradox. The claims presented above would mean that a man who says, "I'm a liar" is neither a truth teller nor a liar, which in the current model is impossible.

However, we must not jump to conclusions. In fact, we do not have paradox here. All that the claims above give is the following proposition:

Proposition: In the model described above (including the added assumption) there can be no person who is able to say the words "I am a liar."

Proof: Let us assume, on the contrary, that there is a man in this model who is able to say "I am a liar." Then, using the claims we have used above, we will arrive at the conclusion that this man cannot be a truth teller and cannot be a liar. The fact that every man in our model must be either a truth teller or a liar will lead us to a contradiction, which in turn leads to the confirmation of the statement above.

The statement we have just proved indicates to us that we may have been hasty in assuming that we were faced with a paradox. Our model (which is indeed hypothetical: in reality there is no such world in which everyone is either a liar or a truth-teller) does indeed contain phenomena that may seem odd or counter-intuitive to some of us. In the strange world we have created there are phrases that are impossible for the people of that world to say. And here we find another important point connected to clear thinking: in models where intuition fails (and there are many such), only the clear formulation of assumptions, definitions and claims will illuminate our path and lead us away from falsehood.

**Achilles and the tortoise:** Imagine a race taking place between Achilles and a tortoise. Since Achilles is faster than the tortoise, he gives the tortoise a head start. Let us track the course of the race according to the following steps. First, Achilles effortlessly reaches point $X_1$, where the tortoise was positioned at the beginning of the race. But in the meantime, the tortoise has reached point $X_2$. When Achilles reaches $X_2$, the tortoise will have moved to point $X_3$, and when Achilles reaches $X_3$, the tortoise will be at $X_4$, etc. For every natural *n*, when Achilles reaches point $X_n$, where the tortoise was seen last, the tortoise will have moved on to point $X_{n+1}$. This process will continue infinitely, and Achilles will never be able to catch up to the tortoise. And yet, we know from experience that this cannot be true. We know that eventually (even if it takes a long time) the faster runner will overtake the slower. We therefore conclude that we are faced with a paradox. What is the root of this paradox?

As we saw in the previous example, if we wish to find the roots of the paradox we must clearly formulate all its assumptions, definitions and claims. First, we should agree on a model in which Achilles and the tortoise are represented by points without size, and the racetrack is represented by a straight line (precisely as things are done in plane geometry). Such a model is necessary if we are to be able to address the location of the two contestants on the track at any given moment. For simplicity's sake, let's assume that the length of the track is a single unit of distance, that Achilles' speed is twice that of the tortoise (which would of course mean either a rather slow Achilles or an exceptionally fast tortoise), and that the tortoise's head start is half the length of the track. Based on these assumptions it is clear to all of us that Achilles should overtake the tortoise precisely at the end of the track, a feat that poor Achilles seems unable to achieve since he cannot traverse the infinite number of points described above (and how could anyone pass through an infinite number of points?).

Now let's move on to examining the process itself. In this process, Achilles traverses half the length of the track (half a unit of length) in the first stage, while the tortoise covers half the distance covered by Achilles (i.e. one quarter of a unit). This means that Achilles will traverse one quarter of a unit in the second stage, and eighth in the third, etc. After *n* stages, Achilles will have traversed a distance of
$\frac{1}{2}+\frac{1}{4}+\frac{1}{8}+\cdots+\frac{1}{2^n}=1-\frac{1}{2^n}$ units of distance (and therefore at no stage will he be able to travel the length of an entire unit). This description confronts us with the fact (or with the assumption that was implicit in the description of the paradox, though it was not stated outright) that the track is divisible into an infinite number of smaller segments (since every natural *n* is assigned a track segment the length of which is $\frac{1}{2^n}$ units). This assumption must also be examined! It may be that such a division is not possible (a not unlikely assumption if we claim that the world is made up of elementary particles that cannot be further divided). If this is indeed the case, then the paradox disappears entirely. The process we have described would end after a finite number of stages, at the last of which Achilles would reach the end of the track and overtake the tortoise (we recommend letting the students present a detailed description of the process).

We are therefore left to discuss what would happen if the track (and time) were indeed infinitely divisible into smaller and smaller segments. If we were to accept such a situation, we would also have to accept that the track itself is made up of an infinite number of points, and that in each stage Achilles traverses (in a finite amount of time) a segment of track that contains an infinite number of points. If we accept this idea, we must also accept that in any finite amount of time Achilles must be passing through an infinite number of points (since we could divide the half of the track Achilles covers in stage one into infinite segments, just as we divided the track as a whole). If we assume that Achilles is able to do so when covering half the track, we must conclude that this ability also applies to the full track. In other words, if we assume that Achilles *can* pass through infinite points in a finite amount of time to complete half the track, he should be able to complete the whole track too: Achilles will complete the whole track in a finite amount of time, and the total distance he will traverse is $\frac{1}{2}+\frac{1}{4}+\frac{1}{8}+\cdots+\frac{1}{2^n}+\cdots=\sum_{k=1}^{\infty}\frac{1}{2^k}=1$ units of distance.

There is no doubt that this conclusion seems strange to anyone unfamiliar with calculus and the theory of limits (i.e. how is it possible to add up an infinite number of terms in a sequence?). But anyone who has studied the theory of limits knows that there is an accurate definition, which is considered satisfactory, for the sum of a sequence with an infinite number of terms (though this may seem odd at first glance). It appears that solving this paradox requires us to define some fundamental concepts. Indeed, approximately 2000 years passed between the time the paradox was first conceived and the time the concept of the limit – which allows us to adequately explain the paradox - was first rigorously defined. It is also worth noting that introducing this paradox while the students are learning about the limit concept will help them appreciate the importance of accurately defining the concept of the limit, and will also help illustrate a concept that is not easily or intuitively perceived: an infinite sequence of numbers which is increasing but bounded.

We conclude this section with another important point. Both of the paradoxes presented above show that clear discussion of the assumptions, concepts and claims is necessary to fully understand the problems. But while in "the Liar" this discussion revealed that the entire paradox was nothing but deception, in "Achilles and the tortoise" it revealed an important, basic issue, namely the need to find a satisfactory definition for concepts like "infinite sum." In doing so, the latter paradox highlights a very basic, difficult and significant problem in the field of mathematics and the exact sciences. In other words, the importance of paradoxes is not limited to their use as "thought exercises"; they can also reveal the existence of basic conceptual problems.

## B: The vagueness problem

**The Heap paradox:** Let us examine the following claims:

1. A collection of 1,000,000 grains of sand is a heap (of sand).

2. If a collection of 1,000,000 grains is a heap, then a collection of 999,999 grains is also a heap.

3. If a collection of 999,999 grains is a heap, then a collection of 999,998 grains is a heap.

Overall, we can make the following claim: "If a collection of a given number of grains of sand constitutes a heap, then when we take away one grain, the remaining collection of grains will also be a heap." This seems to be clear – if we take one grain of sand from a heap, what remains will still be a heap.

The problem begins when we apply the rule many times (for instance one million times on a heap of 1,000,001 grains). This would lead us to the strange conclusion that one grain of sand constitutes a heap. (Worse still, if we apply the rule again we must conclude that the empty set is a heap!)

Have we reached a paradox? Perhaps, though the paradox seems at first glance to be "innocent" and not particularly complex. As we will see in this section and the next, this "innocent" paradox opens the way to an important discussion on the topic of vagueness.

The root of the paradox is in the definition of the concept of "heap". This definition, as perceived in our natural language, is vague. Everyone will agree that a collection of one million grains of sand is a heap, but the average person would not agree that one grain is a heap. And what about ten grains? A hundred? A thousand? It seems clear that in many situations there would be no agreement between different people about whether or not a given collection of grains is a heap. It is this vagueness that complicates our general claim.

Let us conduct a thought experiment in which we take a group of people and place before them a collection of grains that all of them can agree is indeed a heap. Then we will start taking grains away one at a time (it is difficult to imagine anyone having the patience to actually endure the conditions of this experiment, but so long as it is only a thought experiment this need not trouble us). No doubt after a certain number of repetitions, in which everyone will still agree that the collection is a heap, we will arrive at a point in which taking away the next grain leads to disagreement. At least one of the people will claim that what lies before them is no longer a heap, or is no longer definitely a heap, while others may still feel confident that it is. Such a situation must definitely occur so long as everyone is agreed that one grain of sand is not a heap. Thus we will find that there exists a given number of grains *n* that will lead some people to deny (or at least to doubt) the validity of the claim with which we began (since they agreed that *n* grains were a heap, but no longer agreed when another grain was removed).

This situation is characteristic of what happens when a definition is vague. When definitions are vague the limits are not clear: when (i.e. at what number of grains) is a given collection of grains a heap and when is it not? Given the definition's lack of clarity, it is no wonder that certain situations lead to disagreements about its validity.

We must therefore accept that our claim that "if a collection of a given number of grains of sand constitutes a heap, then when we take away one grain the remaining collection of grains will also be a heap" does not rest upon solid foundations. Its current phrasing suggests that it should be valid for any number of grains, but as our experiment showed, this is not the case. Moreover, the claim talks of a "heap" and of the number of grains in a heap without providing us with any knowledge of the limit: at what number of grains does a collection stop being considered a heap? We must therefore realize that it is inadvisable to accept such a claim, the validity of which rests on the validity of a vague definition of "heap." In some situations, (like a million grains) everyone would agree with the claim, while in others they would not. We cannot accept a claim about a heap with a certain number of grains in it when it is unclear to us what number of grains still constitutes a heap. We can generalize and claim here: **any logical claim that relies upon a vague concept could be problematic, since its truth value could depend on the (improperly defined) truth value of the vague concept.**

Let us also take note of what happens if we define a "heap" in the following way (similar to how the concept is defined in computer science):

Definition: Any collection of grains of sand that includes at least one grain is a heap.

Based on this definition, there is no longer any problem with the claim:

"If we have a heap in which there are at least two grains of sand, then if we remove one grain from the heap, the remaining grains will also constitute a heap."

Here the definition of a heap is not vague, the claim is phrased specifically and accurately, and there are no longer any obstacles. On the other hand, when we make a claim based on a vague concept, we are likely to encounter confusion and contradiction.

**The unexpected examination paradox**. A teacher announces to his class that there will be an unexpected examination on one of the days in the following week. The students understand the meaning of an unexpected examination as follows: on the night before the day of the examination, they will not be able to know for sure that the examination will be the next day. Based on that, they make the following series of deductions:

1. The examination cannot be on Friday, since, based on the fact that it must be next week and the fact that it did not take place on previous days, by

> Thursday night everyone would know for sure that the examination would be the next day, which would run counter to the requirements of an unexpected examination.

2. This means that the examination also could not be on Thursday, since, based on the fact that it cannot be on Friday (see previous claim) and the fact that it had not yet occurred, by Wednesday night all the students would know that the examination is the next day. So the examination cannot be on Thursday.

3. The students continue to reason in the same way for the other days, and reach the conclusion that the examination could not take place on the previous days either.

Based on these considerations, the students reach the conclusion that there cannot be an unexpected examination next week, and that the teacher's announcement is fundamentally wrong. On Monday, the teacher arrives in class and announces an (unexpected!) examination. To the students' surprise, they are indeed surprised.

We have a paradox! On the one hand, the students' reasoning seems flawless, so we accept their conclusion that the unexpected examination cannot take place next week. On the other hand, having the unexpected examination that week is clearly possible (since it did indeed take place, surprising everyone).

Before we move on, we should clarify an important point. The teacher told the students two things:

A. There will be an examination next week.

B. The examination will be unexpected (which means, according to the students' interpretation, that on the night before the examination, they will not be able to know for sure that the examination will be the next day).

As we saw in the previous examples, it is important to start by clarifying all of our assumptions. First, we assume that the students have full confidence in the veracity of their teacher's statement (i.e. in part A *and* part B). Without this confidence, the paradox takes on a completely different aspect. Once we have explicitly pointed out this assumption, we can become aware of a problem with the reasoning in item one, namely the conclusion that the night before Friday the students would know for certain that, based on their teacher's statement, the examination would take place the next day. Is this indeed the case? If we only account for part A of the statement and discard part B, this would indeed be our conclusion. But if we fully believe in the statement as a whole, we cannot *at any point* ignore part B! (Especially since part B *is* part of the reasoning later in item 1.) Part B *negates* the possibility of certain knowledge, so we must conclude that, on the night before Friday, full confidence in the (entire) statement leads to an immediate contradiction. In other words, this full confidence is impossible at that point in time! The conclusion described at the end of item 1 above, according to which the examination cannot take place on Friday, was possible only because we initially ignored part B of the teacher's statement and addressed only part A.

While we may come out of this analysis feeling rather confused, this confusion arises from a critical issue in the clear analysis of the paradox – the issue of the students' confidence in both parts of the teacher's statement. We will not expand our direct discussion of this issue here any further (for more on this, and on other aspects of this fascinating paradox, see [KL]). Instead we will move on to focus on the concepts "unexpected" and "certain knowledge" (both of which, of course, are connected to the issue of confidence). Specifically, we will be focusing on the difference between the precise definition of these concepts and their vague definition in natural language (for analysis from other perspectives, see for instance [C], [Sa]).

This paradox seems more difficult than the previous one, and in our opinion this impression is not misleading. Intelligent, experienced people have found themselves perplexed and confused when first faced with this paradox. Indeed, this paradox is quite famous, and dozens of articles have been published about it (see, for example the review in [C]).

A clear analysis does indeed show that, once again, we have fallen victim to our (natural) tendency to use a concept before we have fully clarified its meaning. In this case the "abused" concept is that of the "unexpected examination". In the paradox, the students define (or interpret) the concept of unexpected examination by using the concept of "certain knowledge," but they do not bother to define the latter. Defining one concept by means of another concept, the meaning of which is unclear, cannot be considered as a basis for a full, clear definition. Sticklers may claim (with some justification) that they are not interested in discussing a problem based on concepts that have not been properly clarified. We however, since we like a bit of mystery, will try to examine the paradox anyway. How should we do so?

One way we can try is turning to some sort of probability model (which seems appropriate in light of the reasoning employed by the students). We will illustrate this option now. Let us assume we have a space of five possible events (assuming 5 school days a week), so that for each $k$ ($k$ is a natural number between 1 and 5), the $k^{th}$ event is the existence of the unexpected examination in the $k^{th}$ day. Each such event has a probability of $P_k$ that it will occur (where $P_k$ is a real number that fulfils $0 \leq p_k \leq 1$ as well as $p_1 + p_2 + \cdots + p_5 = 1$). Now let us define the conditional probability $q_k$ for the existence of the $k$ event as the conditional probability we will be calculating - based on additional information gathered or deduced by the evening before day *k* – for the existence or possible existence of the examination on other days. This allows us to determine that the existence of event *k* is a **surprise** (or, equivalently, unexpected) if the calculated value of the conditional probability $q_k$ is smaller than 1 (this definition seems "logical," since we understand that $q_k = 1$ if and only if we have managed to determine with absolute certainty that the examination will take place on the $k^{th}$ day).

This model provides us with a clear definition of the concept "surprise" (and thus also of "unexpected examination"), and also enables us to clearly formulate the students' reasoning.

Indeed, if on Thursday evening the students know that the examination did not take place previously, they will use this knowledge to deduce that the examination must take place on Friday. (While we did point out earlier that there is a fundamental problem with this part of the students' reasoning, based on the assumption that the students believe their teacher's statement in its entirety, we will not, as we have said, be expanding on that here.) They will thus find that $q_5 = 1$, which is impossible given the fact that the existence of the examination – even on Friday – must be a surprise ($q_5 < 1$). The students will thus deduce that the examination could not possibly take place on Friday, and that therefore $p_5 = 0$ (and thus also $p_1 + p_2 + \cdots + p_4 = 1$). If the examination does not take place on days 1-3, the students will deduce, based on this last conclusion ($p_1 + p_2 + \cdots + p_4 = 1$) and on the fact that it did not take place on these three days (which means that $p_1 = p_2 = p_3 = 0$) that on Wednesday evening $q_4 = 1$, so the examination cannot take place on Thursday. The students will use similar reasoning to deduce that the examination cannot take place on the other days either. They could therefore claim that their teacher was wrong: he promised something (an unexpected examination) that cannot exist!

What can we deduce from what we have described above? We have constructed a theoretical model that seems to provide a sufficient definition for the concepts "surprise" and "unexpected examination", and shown that - based on this definition – the unexpected examination cannot take place. Since we know that in reality unexpected examinations *do* take place, we must accept that the model we have constructed, though neat and pleasing, is not connected to reality.

This in turn raises the suspicion that the description of the paradox is making hidden use of two different definitions for the same concept – namely "unexpected examination." One definition is clear and precise, and subject to logical reasoning (like that employed by the students) that can be used to deduce that the unexpected examination cannot occur. But this definition, which may be valid in a theoretical model, does not correspond to the definition of "surprise" or "unexpected examination" in reality. Further thought will lead us to the conclusion that, unlike the concept of "surprise" employed in the model above, the concept of "surprise" used in everyday language is vague, and lacking in accurate definition (as a result, the same is also true of the concepts "certain knowledge" and "unexpected examination"). Moreover, like the concept of the heap, this concept could cause disagreement between different people, since what one person considers a surprise may not seem surprising to another.

Note that at no point during the description of the paradox was any attempt made – by the teacher or the students – to provide a full and accurate definition of the concepts "surprise" and "certain knowledge" (and thus not for "unexpected examination" either). The students' surprise at having the unexpected examination on Monday is therefore their own subjective feeling. Thus, for instance, we (like any of the students) are unafraid of crossing a quiet street at a crosswalk, since in light of our care to follow the rules of the road we are certain of our ability to safely cross the street. If someone were to try and stop us, proving to us that there was still a positive probability that we would be hit by a truck, we would wave off their words with the claim that even if we were to hide in our homes for the rest of our lives (to avoid such dangers) there is still a probability that the ceiling would fall in on our heads.

What this means is that there is a difference between the precise meaning of concepts like "certainty in a given fact" or "certain knowledge" and the meaning we attach to them in daily life. To clarify this issue we will conduct the following thought experiment: the teacher will inform the students that (a) they will have an examination sometime in the next semester and (b) the examination will be unexpected, which means that on the night before the examination the students will not be able to know for sure that the examination will be the next day.

Consider now the following options:

1. The examination has not happened yet and tomorrow is the last day of the semester.

2. The examination has not happened yet and the semester will end in two days.

3. The examination has not happened yet and the semester will end in three days.

4. The examination has not happened yet and the semester will end in four days.

5. The examination has not happened yet and the semester will end in five days.

6. The examination has not happened yet and the semester will end in six days.

7. The examination has not happened yet and the semester will end in ten days.

8. The examination has not happened yet and the semester will end in fifty days.

We ask the students, for each of these options, whether the teacher will have been true to his word (i.e. will have given an "unexpected examination"), assuming the examination takes place the day after the evening in question. There is no doubt that the first option will produce widespread agreement that holding the unexpected examination the next day will not be a surprise, and since the students trust their teacher they will claim that such a situation would be impossible. The teacher is sure to give the examination on one of the earlier days. But what of the other options? The second? The third? The fourth? The fifth?... There is no doubt that at some point one of these other options will elicit a change of opinion amongst some of the students. Some of the students will very likely also be undecided as to whether or not in a given situation the examination can be defined as a surprise. But there is no doubt that all of them would agree that in the eighth situation it would be quite impossible to claim that the examination would not be a surprise. We can now clearly see the similarities between this situation and the Heap paradox, in which it was unclear where we should draw the line at which the definition starts to be valid. What would happen if we asked whether one or two grains were a heap? What if we asked about 10? 100? 1,000? 10,000?

In light of the discussion so far, we find that, like the paradox of the heap, here too the vagueness of the definition is what trips us up. And the definition in this case is that of the concept "surprise" (or alternatively, of "certain knowledge"). The mistake the students made was applying logical considerations, which seem right on their face, to vague concepts. **Logic and vagueness cannot live together in harmony!** (We refer here of course to "classical" logic, which is what the students used, rather than alternatives like "fuzzy logic" and others, which are in themselves noteworthy, but are not in this case what the students used). When we arrive at a vague situation in which the meaning of concepts like "surprise" and "certain knowledge" is unclear, a situation in which one individual can claim certain knowledge and another can claim uncertainty, we must accept that any logical conclusion based on such a state is worthless. That conclusion would waver between "valid" and "invalid" based on the differing meanings attached at any given moment to the concept of "surprise" (or "certain knowledge"). We are therefore faced with the same obstacle that faced us in the "Heap" paradox; the only difference is that here it is hidden and therefore more difficult to identify.

We can sum up this discussion with the following points:

1. We have discovered, to our surprise (surprise? What's that?) that the root of the unexpected examination paradox lies in the vagueness of a concept, which we did not notice in our initial encounter with the paradox. Then, after applying logical reasoning to vague concepts, we were shocked to discover that we have come to conflicting conclusions, without having found any flaw or obstacle in our reasoning.

2. Discussion of this paradox (and many others) teaches us about the importance of **definitions**. If we do not take the trouble to accurately define the concept of "surprise", it can take on different interpretations for different people. Thus one person's understanding of it can lead them to conclude (correctly) that it is impossible for the unexpected examination to take place, while another interpretation (like the one we make use of in natural language) can lead to a different result. It is no coincidence, then, that amongst mathematicians the introduction of every new concept is immediately accompanied by a full and accurate definition.

3. Let us return for a moment to our comment when presenting the paradox of the Heap. In that paradox things seem clear, and discussion seems at first to be unnecessary. On the other hand, the importance of the paradox lies in the clarity with which it introduces the topic of vagueness, which in other cases (like the unexpected examination paradox) is far more hidden and elusive.

4. In light of the previous comment one might ask: what is the point of entering into a discussion of the vagueness of natural language in an article designed to discuss clarity? The answer is clear: although we live perfectly well with vague concepts in real life, such concepts may easily lead us to errors when we try to handle certain problems. Therefore, by discussing vagueness, we can highlight the importance of clear thinking.

Some might say that we have found the root of the unexpected examination paradox and can therefore leave the discussion at that. But have we? It appears that the acceptance of this conclusion depends on the concluder's field of interest (and indeed, we have already noted that this paradox can be addressed from many different points of view). A "zealous" mathematician, who is unwilling to discuss problems that are not clearly phrased, will naturally put an end to the discussion at this point. On the other hand, there are still other avenues to explore that address the logic behind the paradox's story (and these might interest the "zealous" mathematician too). Moreover, many others may feel justified in pursuing the discussion further, since in real life (as we already noted) we manage to deal with vagueness rather well. The truth is that though there is no room for vagueness in mathematics and in the exact sciences, beyond these fields it is nevertheless highly interesting to examine this paradox in the context of natural language. How do we manage so well with the vagueness of the natural language in real life? This philosophical and psychological question is both important and fascinating. The next example illustrates these issues very well.

**The Genie in the Bottle Paradox:** This paradox was based on a famous story by Robert Lois Stevenson [Sh]. The hero of the story is offered the purchase of a bottle in which there resides a genie, who is willing to fulfill almost any wish made by the bottle's owner. On the other hand, the seller of the bottle is faced with a difficult problem: he must sell the bottle at a lower price than that for which he bought the bottle. If he cannot sell the bottle in his lifetime, he will burn in the fires of hell forever (clarification: the genie cannot extend a person's life).

Let us assume, for simplicity's sake, that the price must be set in American dollars. There is no doubt that a rational man would not be willing to buy the bottle for a single cent; nor for two or three. If we continue to employ this line of thinking – like that employed by the students in the unexpected examination paradox – we will arrive at the conclusion that a rational man would be unwilling to buy the bottle for any price (assuming that the buyer is concerned only with the problem of later selling the bottle, and has no other reasons like "I want nothing to do with such dark forces" or "I am willing to risk the fires of hell to save someone I love from a fatal illness").  On the other hand, even if we assume that all of the participants in the story are rational, we would still expect to find a buyer if the price of the bottle was high (e.g. over 1000 dollars). On yet another hand, however, why would anyone find a buyer knowing that one of the buyers (who is also a rational person) would be stuck with the bottle (and with hell) forever? How can the buyer be sure they will be able to sell the bottle? The situation is similar to that of the students in the unexpected examination paradox: can they be certain on Tuesday night that there will be an (unexpected) examination the next day, since the teacher would not "risk" having it on Thursday?

Disagreements about the nature of the paradox, like those described above, arise amongst philosophers as well. Thus, for instance, Quine [Q] notes that the problem in the unexpected examination paradox arises from the vagueness of the concept of "knowledge." He therefore suggests viewing the paradox as a psychological problem rather than a philosophical one. Sorensen [So], on the other hand, claims that in the eyes of epistemologists, "cancelling" the concept of knowledge in this way in order to solve the paradox is like using a bomb to kill a fly. As noted above, we believe that the unexpected examination paradox does present a very interesting epistemological and psychological problem. That said, we would like to add a comment about how this paradox is presented in many different publications. Many analyses of the paradox do not mention the vagueness of concepts like "knowledge" and "surprise" at all, which can be confusing. After all, a clear reading of the wording of the paradox provided us with a "solution" (i.e. **the problem is in the mixing of "classical" logical arguments with vague concepts**), which may lead the reader to wonder – what is the point of going on? To avoid this confusion, we should clearly distinguish between the different analytical approaches.

## C. Interim summary

Let us start by noting some conclusions that should be presented to the students after they have learned about and tried to manage several paradoxes, like those we have described so far. We believe that the analysis of any problem (and especially of paradoxes) must begin with a clear analysis of the problem's concepts and the claims it makes. Therefore, when analyzing a paradox, we must start by clarifying the problem, and by posing questions like:

1. What do the concepts in the paradox mean? Do they have one meaning? Multiple meanings? Are they vague? If they are ambiguous, what is the source of the ambiguity? Can the concepts be interpreted in different ways? What are they? Which of these interpretations are being ascribed to the concepts in this particular problem?

2. What are the claims in this problem? Are they phrased properly? Are they being used properly?

3. Is there a lack of clarity in the problem's presentation? A conflict in the data? A mixture between logical conclusions and intuitive claims? A use of hidden assumptions or information that have not been explicitly stated?

At the end of the process we must arrive at a situation in which we are convinced that we have clarified all of the problem's components. Lack of attention at any point in the process could lead to mistakes, to confusion between the clear and the vague, and to circuitous action that only adds to the confusion. A common mistake, which is not confined to students, is underestimating the complexity and difficulty of the analysis process. Such a process can be complicated and tiring, and even experienced mathematicians can sometimes find that they have failed to clearly understand some component of the problem.

This is therefore a preliminary, but basic and critical, stage in analyzing the paradox. In some cases, this initial analysis can be sufficient to provide a satisfactory explanation of the paradox (as it does, for instance, in the simple version of the Liar paradox above). In other cases, we may find that the clarification process has only revealed a series of new problems (as in the case of Achilles and the Tortoise, where we found that we lacked a method for dealing with problems that involved the concept of infinity).

Another important point we wish to address is the difference between the clear analysis of a paradox (as described in the previous paragraph), and investigation by means of intuitive deduction methods, in which we often make use of vague or ambiguous concepts. Both methods of inquiry have merits of their own, and we believe there is interest to be found in both. Nevertheless, it is important to maintain the distinction between the two, particularly since problems can arise when the second method is applied before the first. In other words, in our opinion any analysis of the paradox must begin with a clear analysis of the type described above. To do otherwise is to risk descending into confusion. When the vague becomes mixed up with the clear, both the former and the latter become vague. This point is discussed in greater detail in the chapter below.

## D. Vagueness vs. Clarity

Timothy Y. Chow [C] opens his article about the unexpected examination paradox with the following words: "Many mathematicians have a dismissive attitude toward paradoxes. This is unfortunate, because many paradoxes are rich in content, having connections with serious mathematical ideas as well as having pedagogical value in teaching elementary logical reasoning." We agree with these words wholeheartedly, and would like to add a further point regarding the dismissive attitude of many mathematicians towards paradoxes. We believe that this attitude arises from two primary sources.

The first source, which we have already noted, is a disinclination to address problems that have confusion and vagueness as their foundations. We disagree with this attitude for the reasons named by Chow. First, paradoxes are not the sole province of confused students. Every thinking human being will encounter (and often fail to avoid) the confusion that arises from the unclear interpretation of a concept. We must look no further than the paradoxes that underlie the very foundations of mathematics (the set theory paradoxes, for instance) for reasons to rethink this attitude. The second reason is that – as we are trying to show in this article – paradoxes have a pedagogical value that should not be discounted.

The second source of mathematicians' dislike of paradoxes is their criticism of how paradoxes are addressed by scholars from other disciplines. Some mathematicians, who value and insist upon clear and rigorous thinking, claim that other disciplines' approach to paradoxes suffers from vagueness and lack of clarity. This suggests that the problem is not with the topic of paradoxes itself, but with how it is handled by those who address it. The literature on paradoxes shows that the claims of such mathematicians are often not wholly unfounded. In many cases we must agree that a "standard analysis" of the paradox is performed before a full analysis of the concepts and assumptions has been completed. Not utilizing this important step to the fullest can impede the clarity and precision of the analysis' later stages. This chapter will be devoted to exploring and illustrating this important point.

This will be achieved through the analysis of two well-known paradoxes: the Raven Paradox and the Grue Paradox. These paradoxes are connected to the concept of "confirmation" of a hypothesis, as described in chapter 5 of [Sa]. Page 92 of the book presents a principle that will serve as the basis for the concept of confirmation throughout our discussion:

G1: A generalization is confirmed by every one of its instances.

Thus, for instance, the generalization "all ravens are black" is confirmed by every instance in which we check a raven and find that it is black. More generally, if a generalization claims that "all A is B," then the presentation, "this A is B" will serve as a confirmation of the generalization. Let us now turn to the first paradox.

**The Raven Paradox.** We will present the paradox, which was introduced by the philosopher of science Carl Hempel in 1945, as it is described in chapter 5 of [Sa].

Let us examine the following principle:

E1: If two hypotheses can be known a-priori to be equivalent (i.e. there is no need for experience to see the equivalence), then any data that confirm one confirm the other.

Here, for instance, are two hypotheses:

R1: All ravens are black.

R2: Everything non-black is non-raven.

These two hypotheses are of course equivalent a-priori. This means that the following instance of R2:

P1: This non-black (in fact, white) thing is non-raven (in fact, a shoe).

Confirms R2, and thus according to E1 also confirms R1, "All ravens are black". This, on the face of it, seems absurd. Data relevant to whether or not all ravens are black must be data about ravens. The color of a shoe can have no bearing whatsoever on the matter. Thus G1 and E1, apparently acceptable principles - lead to the apparently unacceptable conclusion that a white shoe confirms the hypothesis that all ravens are black. This, finally, is our paradox.

The principles of reasoning involved do not appear to be open to challenge, so there are three possible responses:

a) To say that the apparently paradoxical conclusion is, after all, acceptable.

b) To deny E1

c) To deny G1.

This description of the paradox reflects the one provided in [Sa] (which faithfully describes the analyses that have been conducted in important publications in the topic). The book then provides a discussion of the three possible responses. While there is no doubt that this approach to analyzing the paradox is an accepted one, we believe that it skips a step: before we address the popular analysis (i.e. choosing between the three options) we must first clearly examine the basic concepts (especially the concept of confirmation) and the claims that are provided in the paradox. Like the concept of the "heap", this paradox also presents us with the problematic pairing of a vague concept on the one hand, and logical deductions on the other.

We must therefore begin by clarifying the concepts and claims presented in the paradox. Specifically, we will examine:

1. The meaning of the principles and the concepts that were provided.

2. The meaning and validity of the claims from which we deduced that this was indeed a paradox.

Let's start with the latter issue – the claims that led to the conclusion that this is a paradox. These were summed up in the description above as follows: "The color of a shoe can have no bearing whatsoever on the matter. Thus G1 and E1, apparently acceptable principles - lead to the apparently unacceptable conclusion that a white shoe confirms the hypothesis that all ravens are black. This, finally, is our paradox". This statement raises several doubts. The claim that "the color of a shoe can have no bearing whatsoever on the matter" distorts the meaning of the claims that were provided earlier. These claims were designed to confirm the hypothesis "everything non-black is a non-raven" (R2) by examining its instances. Confirmation for this statement is found by looking for objects that are not black and checking whether or not they are ravens. Therefore, such an examination is still relevant to the topic of the raven's color. In other words, it is not the color of a particular shoe that confirms the claim that all ravens are black, but the fact that we have found a white object (or more particularly, a non-black object) and found that it is a shoe (and not a raven).

The confusion here arises, of course, from the fact that in daily life we usually identify the color of a thing at the same time that we identify its nature. The phrasing of R2 requires a certain order in the checking process: if we see a white object, and then examine its nature and find it is a shoe (or more particularly, not a raven), this will count as confirmation for R2. On the other hand, if someone brings us a shoe and suggests checking its color to support R2, we would turn that suggestion down. If we know in advance that the object is a shoe, its color is of no interest to us. Broadly speaking we could say that checking a white object (which turns out to be a shoe and therefore not a raven) would confirm R2, while checking a shoe (no matter what color it turns out to be) would not. The claim "the color of a shoe can have no bearing whatsoever on the matter" therefore distorts the meaning of the examination that is required, and it is only because of this distortion that we conclude the existence of a paradox. Thus we see once again how the vagueness of spoken language ("here is a white shoe!" – which quality did we identify first: an object that is a shoe the color of which must be checked, or a white object the nature of which must be identified?) leads to confusion. Principles and claims must therefore be phrased as clearly and precisely as possible, so as to avoid the fallacies that can arise from linguistic vagueness.

Now let us address the central concept in the paradox – the concept of confirming a hypothesis. This is undoubtedly a vague concept, the meaning of which can give rise to many questions. The evidence one can gather to confirm of a hypothesis by checking a number of instances of that hypothesis in a given group (the group of ravens in hypothesis R1, the group of non-black things in hypothesis R2) can be strong or weak, depending on the size of the group in question. Since the size of the group of non-black things (of any kind) in our world is **immeasurably greater** than the size of the group of ravens, it is not surprising that testing the color of ravens is a much quicker and more efficient method of confirmation than testing the nature of all non-black things – despite the fact that R1 and R2 are a-priori equivalent hypotheses. Moreover, if we were to examine **all** the non-black things in our world and find that all of them are not ravens, we would have gathered full confirmation for the truth of **both** hypotheses. This leads us once again to question whether what we have before us is indeed a paradox, or is rather just a confusion arising from the vagueness of the concept of confirmation.

If we continue thinking clearly about the concept of confirmation, we will find, as many others have found before us, that there are additional problematic aspects to this vague concept, especially in relation to the way it is used in principle G1. We may easily find examples in which "confirmation through instances" does not seem to lend strength to certain hypotheses. We may even find that it is sometimes a source of confusion that can convince us to confirm unreasonable hypotheses. For example, the color of all the parrots in the region where we live is green. We could therefore raise the hypothesis that "all parrots are green." As long as we only make observations in our area and the neighboring towns, we could find a great deal of confirmation for our hypothesis. But if we stop to think for a moment, we will of course reach the conclusion that we should also check the color of parrots in other countries (and we would not even need to go there, a trip to the zoo would suffice) and quickly discover that our hypothesis is incorrect. In other words, if we do not place additional demands and restrictions upon the concept of confirmation, we could spend the rest of our lives finding more and more instances that confirm incorrect hypotheses. Perhaps the central issue here is not the correctness or incorrectness of principles G1 and E1, but the lack of precision in the concept of confirmation that rests at the heart of those principles.

In conclusion, we have seen that if we begin by clearly and patiently examining the concepts and claims laid out in the phrasing of the paradox, and only then attempt to find solutions to it, we will have a clearer picture of the situation, and this will have an impact on the conclusions that we draw. Our preliminary examination of the Raven Paradox, for instance, provided the following conclusions:

1. The claims that led to the conclusion that this was a paradox were faulty.

2. We saw that, at least formally, this was not a paradox.

3. A central issue in the paradox is the vagueness of the concept "confirmation", which was used in principles G1 and E1. It is important to address this *before* we decide if one of these principles is false.

Some of these conclusions are similar to the "accepted" analysis of the paradox as it is presented in [Sa]. However, there is no doubt that the preliminary analysis also revealed new insights, and even presented some of the conclusions from the "accepted" analysis in a different light. We suggest that readers try to reassess the differences by looking again at the detailed description of the paradox in [Sa].

For a sharper, clearer picture of the contribution of conducting preliminary analysis on paradoxes, let us now examine the Grue paradox, which addresses the concept of confirmation. Here too, we suggest that readers compare our analysis to the "accepted" analyses described in [Sa].

**The Grue Paradox (based on [Sa, 5.1.3])**

According to principle G1, green emeralds confirm the hypothesis that all emeralds are green. The Grue Paradox was introduced by Nelson Goodman in 1955 (the word

"Grue" is a combination of the words "green" and "blue"). A given object $x$ is called "Grue" if and only if it meets one of the following conditions:

Gr1: $x$ is green and has been examined, or

Gr2: *x* is blue and has not been examined.

All examined emeralds, being all of them green, count as Grue, according to Gr1. It follows from G1 that the hypothesis that all emeralds are Grue is confirmed by our data: every emerald we have examined is a confirming instance because it was green. This is absurd. If the hypothesis that all emeralds are Grue is true, then unexamined emeralds (supposing that there are any) are blue. This we all believe is false, and certainly not confirmed by our data. G1 must be rejected. What is paradoxical is that a seeming truth, G1, leads by apparently correct reasoning, to a seeming falsehood: that our data concerning emeralds confirm the hypothesis that they are all Grue.

Quite a few articles have been written about this paradox, examining it from different perspectives, the most prominent of which are described in [Sa]. The approaches noted there fall into one of two patterns:

1. Blaming the word "Grue," which is said to be "pathological," and suggesting the need for a general principle to invalidate words like "Grue" so that principle G1 can retain its validity.

2. Claiming that the problem lies not on the word "Grue," but in the attempt to create a principle like G1.

The beginning seems to be good, addressing two central concepts upon which the claim that led us to the absurdity resides. However, this is not enough, in our opinion, for two main reasons:

1. We must begin by going over the entire description of the paradox and examining it clearly. In doing so we will discover, for instance, that aside from the concept "grue" and principle G1 the paradox also describes the **hypothesis** that all emeralds are grue. Maybe the problem is in the hypothesis?

2. Before we begin "placing blame" on a concept or principle, we must first analyze its meaning, its structure etc. We must be careful and not rush to find fault with a concept before we have analyzed it and tried to understand its meaning.

Let us (briefly) examine the phrasing of the paradox. First, the concept "grue" is defined; then the concept is integrated into the hypothesis that "all emeralds are grue."

The concept "grue" was defined thus: object $x$ is grue if it meets the condition, "Gr1 is true **or** Gr2 is true." We should therefore examine Gr1 and Gr2, paying special attention to Gr2, which states that if $x$ is blue and has not been examined, then it is grue. What is the significance of this condition? Does it mean that the examination itself may change the characteristics of $x$ (i.e. make it change color from blue to green)? This would make it a condition that cannot be proved or disproved, and is therefore meaningless. There is after all no way of distinguishing a green emerald that has always been green from a grue emerald that was blue, but has become green because we have observed it!

To further clarify this point, let us raise the following hypothesis:

Hypothesis 1: Every emerald that has not been examined is blue.

We will find that this is a hypothesis that cannot be disproved (and also cannot be proved). The only way to confirm (or deny) such a hypothesis is to look for an emerald that has not been examined and check its color, but the act of checking would make the emerald irrelevant to the hypothesis, since we could claim: how can we know the color of the emerald before we observed it? Maybe its color was different from what we saw? It is no coincidence that the possibility of being disproved is a basic condition required of any scientific hypothesis. In light of all this we can refute the legitimacy of hypothesis 1.

Now let us examine the next hypothesis, which is equivalent to the hypothesis in the paradox (that every emerald is grue).

Hypothesis 2: Every emerald fulfils the condition:

1. It is green and has been examined, or
2. It is blue and has not been examined.

In other words, a trick has been perpetrated here. We have hidden hypothesis 1, which we disqualified, inside another, more complicated claim. Specifically, if we reduce hypothesis 2 so it addresses only unexamined emeralds (and there will of course always be such), hypothesis 2 will produce, either as an instance or as a conclusion, the "illegitimate" hypothesis 1. This means that hypothesis 2 suffers from the same problem as hypothesis 1, except that it is more complex. To clarify this, let us put forth another hypothesis as follows:

Hypothesis 3: The following two hypotheses are valid:

1. Every emerald which has been examined is green

2. Every emerald which has not been examined is blue

We would immediately question the legitimacy of hypothesis 3 in light of the fact that the second condition is invalid (both conditions in hypothesis 3 must be true, and we disqualified condition 2 earlier). A short examination will show that hypothesis 3 is equivalent to hypothesis 2 (and therefore to the hypothesis that every emerald is grue)! So everything that was done with the paradox is a trick: the elusive and confusing word "or" was used to disguise the existence of the invalid hypothesis, hypothesis 1, inside a more complex hypothesis, hypothesis 2.

To conclude: a hypothesis that cannot be either refuted or proved is meaningless. The Grue paradox presents such a meaningless claim, but cleverly disguised through phrasing that uses the confusing connective "or." It is therefore not a paradox we have before us, but a trick. In particular, our "clear" analysis can spare us the need for further analyses based on the two common patterns noted above.

In light of this analysis, and the concerns it raises, one might ask – what is the point of the story about grue? This question could be expanded to include the paradox of the ravens as well, and additional paradoxes that seem to rely largely on the vagueness of concepts and confusion due to phrasing distortion. However, we do not believe that such paradoxes should be discounted completely, for the following reasons:

1. These paradoxes can provide important practice in clear thinking, and not just for students.

2. These paradoxes show us how complex and difficult an issue clear thinking is. In particular, they show us how easy it is to be fooled if we do not properly examine the meaning of important concepts and claims.

3. The Grue paradox shows that the use of logic is not a natural, automatic thing. Even the use of the innocent connective "or" can be a stumbling block (would a statement like "it is always true that 3 is less or equal to 3 be immediately accepted as valid by anyone?). This is an important point to internalize in teaching: students will often have difficulty managing the formal side of logic in general, and particularly connectives and quantifiers.

4. Discussing and debating these last two paradoxes (Ravens and Grue) could contribute a great deal to the study of philosophy of science.

5. These paradoxes, and others like them, can teach us about our intuitive way of thinking, about the problems associated with the thinking we call

"rational," about how we perceive concepts, etc. Clear thinking analysis can show us how easy it is to make mistakes, or to "miss the point," if we do not pay proper attention and think about the problem clearly.

Finally, describing and analyzing paradoxes shows us that the clear analysis of problems can be a challenge for experts, for philosophers and, as we will see in the next chapter, for famous mathematicians as well.

## E. The clear thinking problem: not just for students

It is not easy to convince students of the importance of clear thinking. Moreover, their willingness to learn and practice is further hampered by the fact that they often have the impression that this is a relatively simple skill to acquire – that all we must do is learn not to be hasty, to carefully check phrasing, to look again...and that is all. Unfortunately, things are not really that simple. We therefore think it is important to present students with examples that show how even the most well-known and respected experts can find themselves struggling with issues concerning the clarity of how a problem or a definition is presented.

A typical example of this is Russel's paradox. For many years, mathematicians did not bother to look for a rigorous definition of the basic concept of a set. This may seem odd – after all, a set is a basic concept in mathematics! On the other hand, however, it may not be so surprising. This is a concept that seems very intuitive and easy to grasp, suggesting that it can be accepted as a basic concept without the need for an accurate definition (like the concept of a point in plane geometry). What might still seem strange is the fact that even after paradoxes associated with this concept started to arise (like Cantor's paradox), no real efforts were made to place the concept of a set upon firm foundations. It seems that the mathematicians of the time were satisfied with the intuitive impression that there "isn't really a problem" after all with the concept. Thus, for instance, Gottlob Frege, one of the fathers of modern logic, began a far reaching and lengthy project by establishing mathematics on its most basic foundations, of which the concept of a set was a central one. He had published the first volume of this work, and was about to publish the second, when Bertrand Russel released his famous paradox, which we will describe briefly below.

A set can include another set as an element. More specifically, a set can include itself as an element (e.g. the set of all sets that can be described in less than 100 words is an element in itself because it can be described in less than 100 words). Let us define a set *A* as the set of all sets that do not include themselves as an element, and ask if *A* includes itself as an element. On the one hand, if *A* includes itself as an element, then the definition of *A* will tell us that it does not include itself as an element, a contradiction. On the other hand, if A does not include itself as an element, then the definition of A will tell us that it does include itself, and we have a contradiction once again. And thus we have a paradox: we have created a set that cannot include itself as an element, but also cannot exclude itself.

This simple paradox undermines the concept of a set, which is a cornerstone in the construction of mathematics as a whole! Indeed, the publication of the paradox led Frege to abandon his life's work. The problem raised by the paradox demanded a solution. And so, after many years of effort from some of the world's leading mathematicians, who attempted a variety of solutions, we seem to have found an acceptable definition for the concept of a set. It is a complex, constructive definition, which is usually taught in advanced mathematics courses.

The example we have provided here can be embellished with additional pertinent examples from the history of mathematics, describing the bumpy road mathematics had to travel before arriving at the rigorous forms and descriptions that exist today. But there is no need to delve too far into the history of mathematics. Every mathematician can recall many times when he or she was forced to give up on an idea of which they had been certain, when the attempt to write out the necessary definitions and proofs in accurate detail revealed previously unnoticed problems. Only through full and accurate documentation can we monitor the consistency and the validity of our work. There have been many cases in which errors were found in mathematical articles that had gone through a full process of review – both by the authors of the article and by their peers. It is no wonder that often these errors are found in articles written by authors who neglected to provide full written details, writing instead that "one can easily see that…" We could also add stories (and there is no lack of these) about false proofs that have been published by expert mathematicians.

We hope that the above may help convince students to acknowledge the importance, as well as the difficulty, of adopting the habits of clear thinking. But that is not enough. We believe that this acknowledgement is being overlooked by the educational system – from elementary school to the university – and that as a result insufficient effort is being put into providing students with such skills and habits while they are in the system. A lengthy discussion of this extremely important problem is, unfortunately, beyond the scope of this paper (for more details, see [L], and also [EW]).

## F. The clear thinking problem: not just for mathematicians

Acknowledging the importance of clear thinking, and especially acknowledging the difficulty of acquiring clear thinking skills, is not just important in mathematics, but in other fields as well. After all, clear understanding of text, understanding the meaning of concepts, being able to state a claim simply and correctly, are all important in any field, from the exact sciences to the humanities. One often hears complaints from lecturers on the humanities and social sciences that the students have a tendency to "jump the gun," trying to learn complex theories and without having first acquired the necessary foundations, without having gained proficiency in properly reading and understanding a text. University lecturers also complain about the level of knowledge and skill with which their students arrive at the university. These complaints do not refer to the students' ability to understand and implement complex theories, but once again, they lament their inability to understand a text, to comprehend the significance of concepts, to write clearly – in short, the students' lack of the skills associated here with the term "clear thinking."

Most universities are also not making a proper effort to remedy the situation (doing so would not be simple; it requires a great investment of time and resources). There seems to be some sort of unjustified "faith" that students will be able to acquire the needed capabilities independently.

The issue is also relevant to the topic of academic research and publication. The field of mathematics employs a method of work and review that has been developed over many years. Every concept in mathematics must be defined clearly and unequivocally; every claim must be properly phrased, and will not be accepted without proofs that are stated in a manner that meets accepted criteria. This method of review, which is widely accepted and extremely strict, does not exist in many other academic fields. This is unsurprising, since while all concepts in mathematics are unambiguous and specific, and its laws have a clear truth value, the same cannot be said of many other disciplines. Broadly speaking we could say that the further we get from mathematics, through the hard sciences, towards the social sciences and then the humanities, we will meet with more and more concepts that are vague or ambiguous, and with claims that have no definite truth value.

One might conclude from this that such inexact realms as those inhabited by the humanities would therefore have no need of the "clear" methodologies of mathematics, but we wish to dispute that claim. In areas where it is easy to become confused by the various meanings of a concept and to assign an unequivocal truth value to a vague claim, it is *especially* important to maintain a strict method of critical appraisal. The paradoxes we examined in this article are just a taste of how this critical method can be applied, how it can prevent the confusion that arises from vague terminology and convoluted, misleading phrasing. This method is not always incorporated in the "accepted approaches" to solving paradoxes like these. One might therefore ask whether it would not be proper for those who work on paradoxes to borrow some of the mathematicians' methods and apply them to the early stages of their problem analysis. We believe this question is relevant to many additional fields as well, including the humanities. It is in the humanities where we often see cases in which the style and impact of the arguments is emphasized, at the expense of their validity and coherence. Moreover, these fields are the site of bitter disagreements between experts. One might therefore ask, would it not be possible to try and create some basic critical system, which would at least provide boundaries that are accepted by all those who work in the humanities – or at least most of them? Creating such a system is of course a complex and difficult task, but we believe, in light of the current state of the humanities and the extreme disagreements that exist there, that such an effort would be worthwhile. Such a unified and accepted critical system would reinforce existing knowledge, clarify phrasing and monitor the validity of the use of concepts and claims. We believe it is important to try and assess this possibility in order to strengthen, if only slightly, the boundaries that allow us to distinguish between claims that have worth and relevance, and claims that do not.

**Arieh Lev** School of Computer Sciences, The Academic College of Tel-Aviv-Yafo, Tel-Aviv, Israel. arieh@mta.ac.il

**Gil Kaplan** School of Computer Sciences, The Academic College of Tel-Aviv-Yafo, Tel-Aviv, Israel. gilk@mta.ac.il